\documentclass[3p,times]{elsarticle}

\usepackage{ecrc}

\volume{00}
\firstpage{1}

\journalname{Journal}

\runauth{}

\jid{Journal}

\jnltitlelogo{Journal}

\usepackage{amssymb}

\usepackage[figuresright]{rotating}
\usepackage{dsfont}
\usepackage{tabu}
\usepackage{amsmath}
\usepackage{amsthm}
\usepackage{dsfont}
\usepackage{amsfonts}
\usepackage[subnum]{cases}
\usepackage{amsmath,graphicx}
\usepackage{amssymb}
\usepackage{graphicx}
\usepackage{epstopdf}
\usepackage{tabularx}
\usepackage{color}
\usepackage{times}
\usepackage{epsfig}
\usepackage{amsmath}
\usepackage{array}
\usepackage{threeparttable}
\usepackage{algorithm}
\usepackage{color}
\usepackage{lineno}
\usepackage{bm}
\usepackage{subfigure}
\usepackage{cases}
\usepackage{epsfig,amssymb,subfigure,amsmath,version}
\usepackage{amssymb,version,graphicx,fancybox,mathrsfs}
\usepackage{graphicx,mathrsfs,fancyhdr}
\usepackage{subfigure,slashbox,rotating,version,fancybox,pifont,booktabs,epsfig}
\newtheorem{theorem}{Theorem}[section]
\newtheorem{remark}{Remark}[section]

\begin{document}

\begin{frontmatter}



\dochead{}

\title {Primal-dual method to the minimized surface regularization  for image restoration}


\author{Zhi-Feng Pang$^\dag$ and Yuping Duan$^\ddag$}

\address{$^\dag$ College of Mathematics
and Statistics, Henan University, Kaifeng, 475004, China.}
\address{$^\ddag$  Center for Applied Mathematics, Tianjing  University, Tianjin, 300072, China. }
\ead{zhifengpang@163.com and doveduan@gmail.com}
\begin{abstract}
We propose a new image restoration model based on the minimized surface regularization. The proposed model closely relates to the classical smoothing ROF model \cite{4}. We can reformulate the proposed model as a min-max problem and solve it using the primal dual method. Relying on the convex conjugate, the convergence of the algorithm is provided as well. Numerical implementations mainly emphasize the effectiveness of the proposed method by comparing it to other well-known methods in terms of the CPU time and restored quality.
\end{abstract}


\end{frontmatter}



\section{Introduction}
Image restoration  is one of the most fundamental and important problems in low-level image processing, which is the operation to recover (as good as possible) the clean image $u:\Omega\rightarrow \mathbb{R}$, $\Omega\subset\mathbb{R}^2$, from a contaminated image $f:\Omega\rightarrow \mathbb{R}$ as
\begin{eqnarray*}
f=Ku+\eta,
\end{eqnarray*}
where $K$ is linear degraded operator (blur operator) and $\eta$ is an additive noise. An ideal restored model is expected to enhance image by reducing degradations in the image acquisition process and preserving edges as much as possible. However, it is often difficult to simultaneously remove the noise and enhance edges because both the noise and edges are high frequency signals.

During the past several decades, the models based on variational partial differential equation  (PDE) have been attracted much attention such as TV-based models \cite{4, 5,10} and nonlocal-based models \cite{26,27} and also obtained some satisfactory results \cite{1,2,3}. Different to  aforementioned models, which are developed based on the image domain, the authors in \cite{28,29,38} proposed to consider the image as an embedded surface $\mathcal{M}\in\mathbb{R}^3$ denoted by
\begin{eqnarray*}
\Omega\rightarrow \mathcal{M}:\mathrm{x} \rightarrow\mathcal{U}(\mathrm{x}),
\end{eqnarray*}
where $\mathrm{x}:=(x_1,x_2)$ denotes the local coordinates of the surface and $\mathcal{U}(x_1,x_2):=\big(x_1,x_2,u(x_1,x_2)\big)$. Note that $\Omega$ and $\mathcal{M}$ are viewed as Riemannian manifold equipped with suitable metrics. By introducing metrics $\mathrm{d}^2s\alpha\mathrm{d}^2x_1+\alpha\mathrm{d}^2x_2$ on $\Omega$ and $\mathrm{d}^2\tilde{s}=\alpha\mathrm{d}^2x_1+\alpha\mathrm{d}^2x_1+\mathrm{d}^2u$ on $\mathcal{M}$,
we can obtain
\begin{eqnarray}\label{11}
\mathrm{d^2\tilde{s}}=(\mathrm{d}x_1,\mathrm{d}x_2)
\left[\begin{array}{cc}\alpha+u_{x_1}^2 & u_{x_1}u_{x_2}\\u_{x_1}u_{x_2}&\alpha+u_{x_2}^2
 \end{array}\right]
\left[\begin{array}{cc}\mathrm{d}x_1\\ \mathrm{d}x_2
 \end{array}\right],
\end{eqnarray}
where $\alpha>0$ is a shrinkaging parameter for the local coordinates $(\mathrm{d}x_1,\mathrm{d}x_2)$ and $\mathrm{d^2(\cdot)}$ denotes $\big(\mathrm{d(\cdot)}\big)^2$ with the convention. In order to obtain a restored approximation $u$ from $f$, we need to search for $\mathcal{M}$ with the minimal area. In this way, singularities are smoothed. Let $g$ denote the determinant of the second-order square matrix in (\ref{11}). We consider to minimize $u$ as follows
\begin{eqnarray}\label{12}
\min_{u}~~\mathcal{J}_{\alpha}(\nabla u):=\int_{\Omega}\sqrt{g}\mathrm{dx}=\sqrt{\alpha}\int_{\Omega}\sqrt{\alpha+|\nabla u|^2}\mathrm{dx}.
\end{eqnarray}

The Euler-Lagrange equation of (\ref{12}) gives
\begin{eqnarray}\label{13}
E_{\alpha}(u)=0,
\end{eqnarray}
where $$E_{\alpha}(u):=\mathrm{div}\left(\frac{\nabla u}{\sqrt{|\nabla
u|^2+\alpha}}\right).$$

It is clear that the mean curvature of $\mathcal{M}$ is zero when $\alpha=1$. Surfaces of zero mean curvature are known as minimal surfaces. Thus, we can solve (\ref{13}) by embedding it into the following dynamical scheme
\begin{eqnarray*}
\frac{\mathrm{d\mathcal{X}}}{\mathrm{dt}}(t)=E_{\alpha}(u),
\end{eqnarray*}
where $\mathcal{X}(t)=(x_1,x_2,u(t,x_1,x_2))$. However, this scheme only considers how to regularize the image while ignores to preserve the image features. Therefore, we propose a novel model by introducing a data fitting term as follows
\begin{eqnarray}\label{14}
\min_{u}~~\frac{\lambda}{2}\left\|Ku-f\right\|_2^2+\int_{\Omega}\sqrt{\alpha+|\nabla u|^2}\mathrm{dx},
\end{eqnarray}
where $\lambda>0$ is a positive parameter and $\|\cdot\|_2$ denotes the $L^2$-norm.

The proposed model (\ref{14}) closely relates to the classic ROF model proposed by Rudin, Osher and Fatemi (ROF model) \cite{4} when $\alpha=0$.  On the other hand, when $\alpha>0$ as in our model (\ref{14}), extra smoothness is introduced to the Total Variation (TV). We can employ similar numerical methods to solve the smoothing ROF method (\ref{14}) as the classical ROF model including the time marching scheme \cite{4} and the fixed point iteration scheme \cite{8}. These methods are usually restricted to the Courant-Friedrichs-Lewy (CFL) condition and the data scale of the operator inversion.  In this paper, we propose a primal-dual method to solve the model (\ref{14}). To the best of our knowledge, this method has not been used to solve the model (\ref{14}) although it has been verified on some non-smoothing models in the field of image processing and machine learning. Moreover, we firstly employ the Legendre-Fenchel transformation to reformulate the minimization problem (\ref{14}) as a saddle-point problem and use the alternative updating scheme to solve the primal and dual variables. The proposed primal-dual algorithm can help to avoid the difficulties when working solely with the primal variable or dual variable \cite{4,7}. We use numerical experiments to demonstrate the proposed primal-dual algorithm can achieve the solution in a reasonable time.

The contents of the paper are arranged as follows.  In section 2, we give some preliminaries of the primal-dual method and use it to solve the proposed model. Some numerical comparisons are done between the proposed method and other classic numerical methods in section 3. We give the concluding remarks in section 4.

\section{The basic results}
Firstly, we define the convex conjugation as \cite{15} in the following way
\begin{eqnarray}\label{21}
H^*(s)=\sup_{t}\left\{\langle s,t\rangle-H(t)\right\}
\end{eqnarray}
for a function $H:\mathcal{Z}\rightarrow\mathds{R}\cup\{\infty\}$ in order to obtain the saddle-point problem from the model (\ref{14}).  Here, $\mathcal{Z}$ denotes a finite-dimensional Hilbert space. Throughout the rest of the paper, we assume the images are matrices with the size of $N\times N$ and with the periodic boundary condition.  Let us define the Euclidean space $X=R^{N\times N}$ and $Y=X\times X$. The usual scalar products can be denoted as $\langle \textbf{v},\textbf{w}\rangle_X:=\displaystyle\sum^n_{i=1}\sum^n_{j=1}\textbf{v}_{i,j}\textbf{w}_{i,j}$ with the  norm $\|\textbf{v}\|_X=\sqrt{\langle \textbf{v},\textbf{v}\rangle_X}$ for $\textbf{v},\textbf{w}\in X$ and $\langle \textbf{p},\textbf{q}\rangle_Y:=\displaystyle\sum^n_{i=1}\sum^n_{j=1}\sum^2_{\iota=1}\textbf{p}^{\iota}_{i,j}
\textbf{w}^{\iota}_{i,j}$ for $\textbf{p},\textbf{q}\in Y$.  In the following, we define the discrete gradient $\nabla u=(D^+_xu,D^+_yu)$ with the forward difference operators
\begin{eqnarray*}
D^+_xu_{i,j}=\left\{
\begin{array}{ll}
u_{i+1,j}-u_{i,j}\hspace{3pt}&\mbox{if}\hspace{3pt}1\leq i<n, 1\leq j\leq n,\\
u_{1,j}-u_{i,j}\hspace{3pt}&\mbox{if}\hspace{3pt}i=n, 1\leq j\leq n,
\end{array}
\right.\\
D^+_yu_{i,j}=\left\{
\begin{array}{ll}
u_{i,j+1}-u_{i,j}\hspace{3pt}&\mbox{if}\hspace{3pt}1\leq i\leq  n, 1\leq j<n,\\
u_{i,1}-u_{i,j}\hspace{3pt}&\mbox{if}\hspace{3pt}1\leq i\leq  n, j=n.
\end{array}
\right.
\end{eqnarray*}
We also define the backward difference operators as
\begin{eqnarray*}
D^-_xp^1_{i,j}=\left\{
\begin{array}{ll}
p^1_{i,j}-p^1_{i-1,j}\hspace{3pt}&\mbox{if}\hspace{3pt}1<i\leq n, 1\leq j\leq n,\\
p^1_{i,j}-p^1_{n,j}\hspace{3pt}&\mbox{if}\hspace{3pt}i=1, 1\leq j\leq n,
\end{array}
\right.\\
D^-_yp^2_{i,j}=\left\{
\begin{array}{ll}
p^2_{i,j}-p^2_{i,j-1}\hspace{3pt}&\mbox{if}\hspace{3pt}1\leq i\leq  n, 1<j\leq n,\\
p^2_{i,j}-p^2_{i,n}\hspace{3pt}&\mbox{if}\hspace{3pt}1\leq i\leq  n, j=1.
\end{array}
\right.
\end{eqnarray*}
Based on the relation $\langle u, \mathrm{div}\textbf{p}\rangle_X=-\langle\nabla u, \textbf{p}\rangle_Y$ in \cite{6}, we can obtain the divergence operator as
\begin{eqnarray*}
\mbox{div}\textbf{p}=\nabla_x^-p^1+\nabla_y^-p^2
\end{eqnarray*}
for $\textbf{p}\in Y$.
Therefore, the discrete equivalent of (\ref{14}) can be denoted by
\begin{eqnarray}\label{22}
\min_u\left\|Ku-f\right\|_X^2+\left\|\sqrt{\alpha+|\nabla u|^2}\right\|_X,
\end{eqnarray}
where $|\nabla u|^2=\left(D_x^+u\right)^2+\left(D_y^+u\right)^2$. Following from the convex conjugate (\ref{21}) (See Exam 8.5 in  \cite{11}), there is
\begin{eqnarray*}
g(t)=\sqrt{\alpha+t^2}\Longleftrightarrow g(t)=\sup_{|s|\leq1}\left\{\langle t, s\rangle+\sqrt{\alpha\left(1-s^2\right)}\right\}.
\end{eqnarray*}
Thus, we can rewrite the minimization problem (\ref{22}) as a min-max problem
\begin{eqnarray}\label{23}
\min_u\max_{\|\textbf{p}\|_{\infty}\leq1}~\left\{\frac{\lambda}{2}\|Ku-f\|_X^2-\langle u,\mathrm{div}\textbf{p}\rangle_X+\left\|\sqrt{\alpha\left(1-|\textbf{p}|^2\right)}\right\|_X\right\},
\end{eqnarray}
for $\textbf{p}\in Y$.
Since the subjective function (\ref{23}) is proper convex, we can interchange the order of min and max and solve the problem using the primal-dual scheme \cite{7}. We separate (\ref{23}) into the following two subproblems.
\begin{itemize}
\item[$\bullet$]For the primal variable $u$  in (\ref{23}): By ignoring the unrelated term to $u$, we can obtain
\begin{eqnarray*}
\min_u~~\left\{\frac{\lambda}{2}\|Ku-f\|_X^2-\langle u,\mathrm{div}\textbf{p}\rangle_X\right\},
\end{eqnarray*}
with its optimality condition as
\begin{eqnarray*}
\lambda K^T(Ku-f)-\mathrm{div}\textbf{p}=0,
\end{eqnarray*}
where $T$ denotes the matrix transpose.
In general, the blurring operator matrix $K$ is ill-posed,  we can use the gradient method to compute $u$
\begin{eqnarray*}
u^k-u^{k+1}=\tau\left(\lambda K^T(Ku^{k+1}-f)-\mathrm{div}\textbf{p}\right).
\end{eqnarray*}
Due to the assumption of the periodic boundary condition, we can use the Fast Fourier Transform (FFT) to solve $u$ as
\begin{eqnarray}\label{24}
u^{k+1}_{\textbf{p}}:=u^{k+1}=\mathcal{F}^{-1}\left(\frac{\mathcal{F}(u^k)+\lambda\tau \mathcal{F}(K^Tf)+\tau \mathcal{F}(\mathrm{div}\textbf{p})}{\mathcal{F}(I)+\lambda\tau\mathcal{F}(K^TK)}\right),
\end{eqnarray}
where $\mathcal{F}^{-1}$ denotes the inverse transform of $\mathcal{F}$ and $I$ is the identity matrix.
\item[$\bullet$]For the dual variable $\textbf{p}$ in (\ref{23}): By ignoring the unrelated term to $\textbf{p}$ and introducing an indicator function
\begin{eqnarray*}
\chi_{\mathcal{K}}\left(\textbf{p}\right)=\left\{
\begin{array}{ll}
0,&\mbox{if}\hspace{5pt}\textbf{p}\in\mathcal{K},\\
+\infty,&\mbox{if}\hspace{5pt}\textbf{p}\not\in\mathcal{K},
\end{array}\right.
\end{eqnarray*}
where $\mathcal{K}:=\left\{|p|_{\infty}\leq1\right\}$.  Then we have
\begin{eqnarray*}
\max_{\textbf{p}}~~\left\{\langle\nabla u,\textbf{p}\rangle_Y
+\left\|\sqrt{\alpha\left(1-|\textbf{p}|^2\right)}\right\|_X+\chi_\mathcal{K}(\textbf{p})\right\}.
\end{eqnarray*}
The optimality condition of the above maximization problem is
\begin{eqnarray*}
\left(\nabla u+\partial\chi_\mathcal{K}(\textbf{p})\right)\sqrt{1-|\textbf{p}|^2}-\sqrt{\alpha}\textbf{p}=0.
\end{eqnarray*}
Note that $\partial\chi_\mathcal{K}(\textbf{p})=0$ because indicator function is a constant function.  Here, $\partial$ denotes the sub-gradient defined by $\partial\hbar(\bar{y}):=\{v|\hbar(\bar{x})-\hbar(y)\geq(v,\bar{x}-\bar{y})\}$ at the point $\bar{y}$ for a function $\hbar$. Therefore, using the projection gradient method, we have
\begin{eqnarray}\label{25}
\textbf{p}^{k+1}_u:=\textbf{p}^{k+1}=\frac{\textbf{p}^k+\sigma\left(\nabla u
\sqrt{1-\left|\textbf{p}^k\right|^2}-\sqrt{\alpha}\textbf{p}^k\right)}{\max\left\{1,~
\left|\textbf{p}^k+\sigma\left(\nabla u\sqrt{1-\left|\textbf{p}^k\right|^2}-\sqrt{\alpha}\textbf{p}^k\right)\right|\right\}}.
\end{eqnarray}
\end{itemize}

\begin{theorem}\label{th21}
Assume that $\tau\sigma<1/8$ and choosing $u^1=f$ and $\overline{\textbf{p}}^k=\textbf{0}$,
then the sequence $\left\{\left(u^{k+1},\textbf{p}^{k+1}\right)\right\}$ generated by \begin{eqnarray}\label{26}
\left\{
\begin{array}{ll}
u^{k+1}=\mathcal{F}^{-1}\left(\frac{\mathcal{F}(u^k)+\lambda\tau \mathcal{F}(K^Tf)+\tau \mathcal{F}(\mathrm{div}\overline{\textbf{p}}^k)}{\mathcal{F}(I)+\lambda\tau\mathcal{F}(K^TK)}\right);\\
\overline{u}^k \hspace{7pt}= 2u^{k+1}-u^k;\\
\textbf{p}^{k+1}=\frac{\textbf{p}^k+\sigma\left(\nabla\overline{u}^k
\sqrt{1-\left|\textbf{p}^k\right|^2}-\sqrt{\alpha}\textbf{p}^k\right)}{\max\left\{1,~
\left|\textbf{p}^k+\sigma\left(\nabla \overline{u}^k\sqrt{1-\left|\textbf{p}^k\right|^2}-\sqrt{\alpha}\textbf{p}^k\right)\right|\right\}};\\
\overline{\textbf{p}}^k\hspace{5pt}=\textbf{p}^{k+1};
\end{array}\right.
\end{eqnarray}
converges to the saddle point $\left(u^*,\textbf{p}^*\right)$ of the problem (\ref{23}). Furthermore, $u^*$ is the solution of the  problem (\ref{22}).
\end{theorem}
\begin{remark}
For Theorem \ref{th21}, if we set $F(u)=\frac{\lambda}{2}\|Ku-f\|_X^2$ and $G(\nabla u)=\left\|\sqrt{|\nabla u|^2+\alpha}\right\|_X$, the iteration scheme (\ref{26}) is the exact Algorithm 1 used in \cite{7}. Furthermore, the convergence can be kept since $\|\nabla\|_Y\leq 1/8$ (See Theorem 1 in \cite{7}). Note that the operator norm $\|\cdot\|_Y$ is defined as $\|A\|_Y=\max \big\{\|A\textbf{z}\|_Y,\hspace{3pt}\textbf{z}\in Y~ \mbox{with}~ \|\textbf{z}\|_Y\leq1 \big\}$.
\end{remark}

\section{Numerical implementations}
In numerical implementations, we consider to use the proposed model (\ref{14}) for the basic image restoration problems, i.e., image denoising and deblurring. In fact, the model (\ref{14}) has many other applications, for example the CT or MRI medical image reconstruction problems and image inpainting problems with different operators $K$, etc.. In order to demonstrate the advantage of the proposed Primal-Dual Method (PDM), we compare it with another two classic numerical methods, which are the Time Marching Method (TMM) \cite{4} and the Fixed Point Method (FPM) \cite{5}. All the algorithms will stop when  $\max\left\{\frac{\left\|u^{k+1}-u^{k}\right\|_X}{\left\|u^k\right\|_X}, \frac{\left|E(u^{k+1})-E(u^{k})\right|}{\left|E(u^k)\right|}\right\}\leq10^{-5}$ or the iteration arrives to 500. The simulations are preformed in Matlab 7.14(R20014a) on a PC with an Intel Core i5 M520 at 2.40 GHz and 4 GB of memory.

\begin{figure*}[ht!]
\centering
\begin{minipage}[htbp]{0.18\linewidth}
\centering
\includegraphics[width=1.2in]{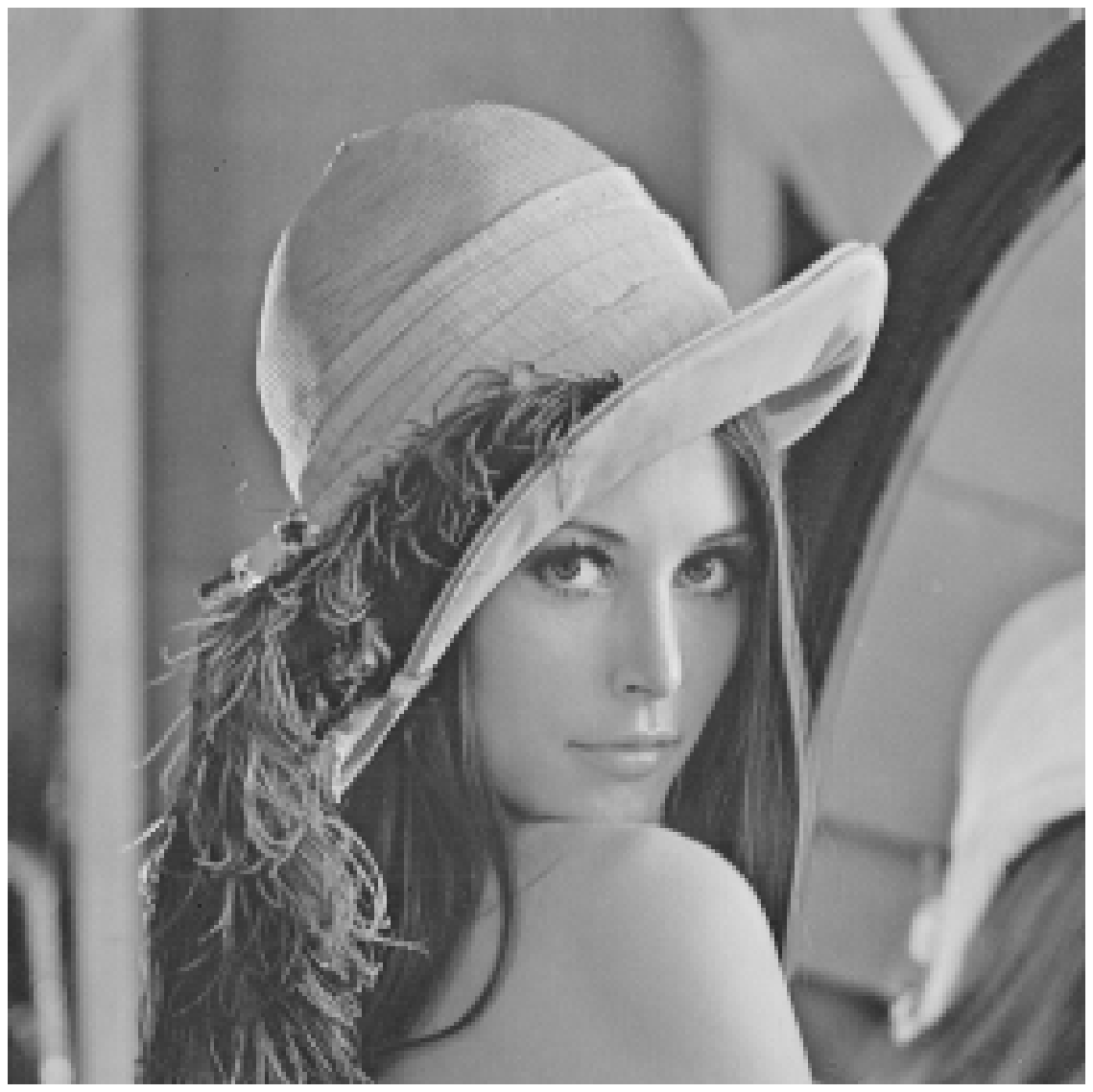}
{(a) Lena}
\end{minipage}~
\begin{minipage}[htbp]{0.18\linewidth}
\centering
\includegraphics[width=1.2in]{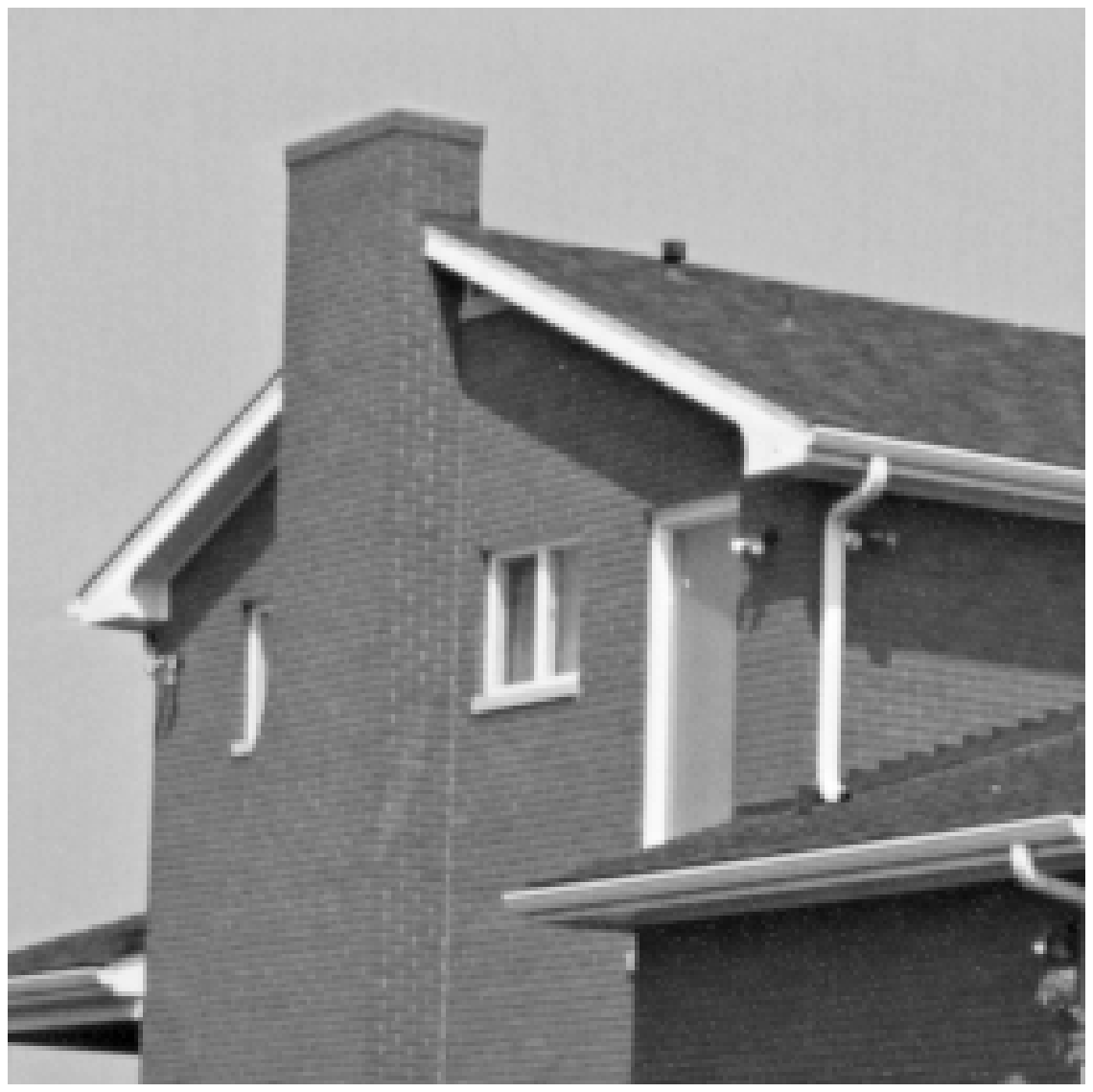}
{(b) House}
\end{minipage}~
\begin{minipage}[htbp]{0.18\linewidth}
\centering
\includegraphics[width=1.2in]{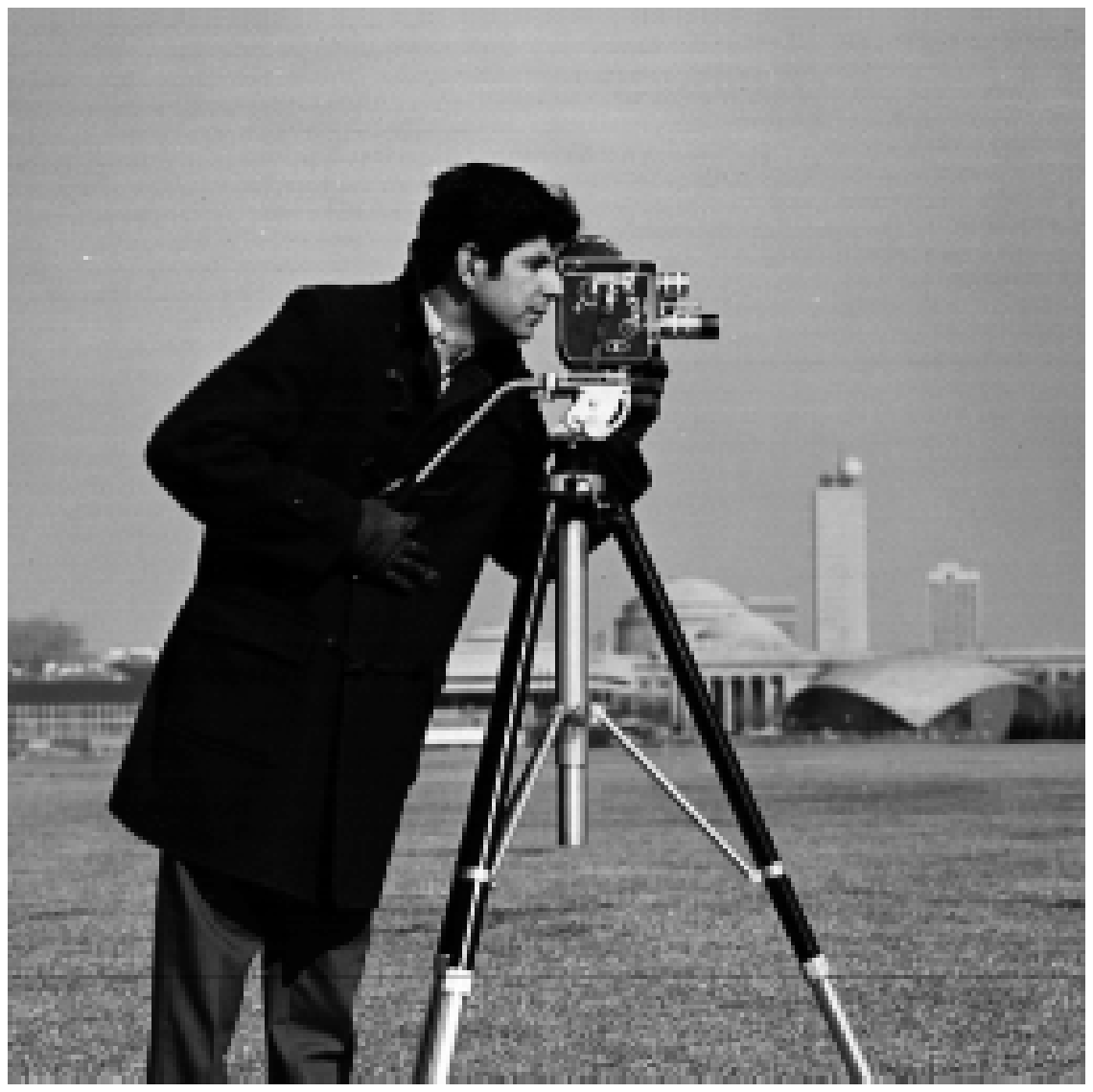}
{(c) Cameraman}
\end{minipage}~
\begin{minipage}[htbp]{0.18\linewidth}
\centering
\includegraphics[width=1.2in]{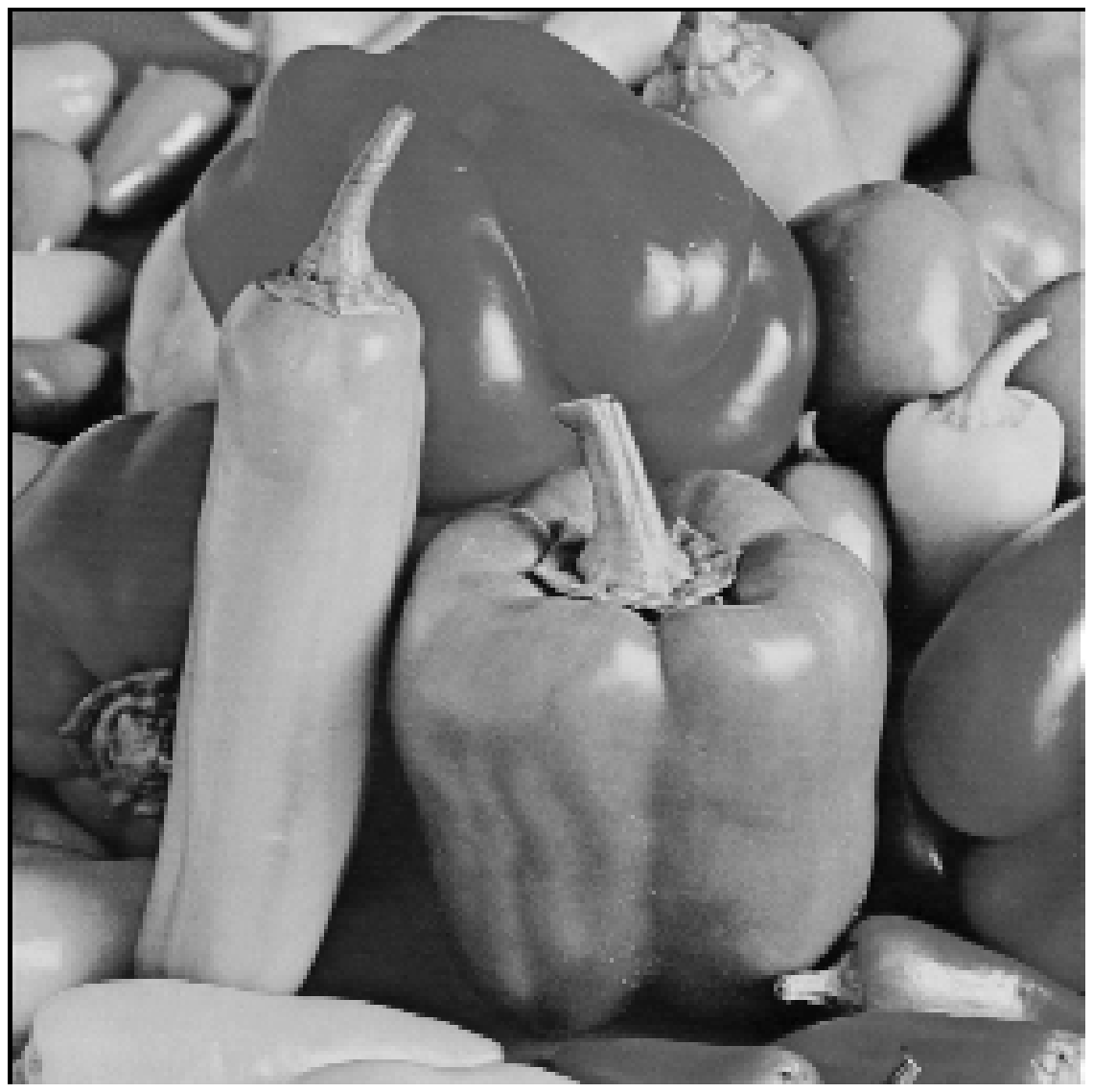}
{(d) Peppers}
\end{minipage}~
\begin{minipage}[htbp]{0.18\linewidth}
\centering
\includegraphics[width=1.2in]{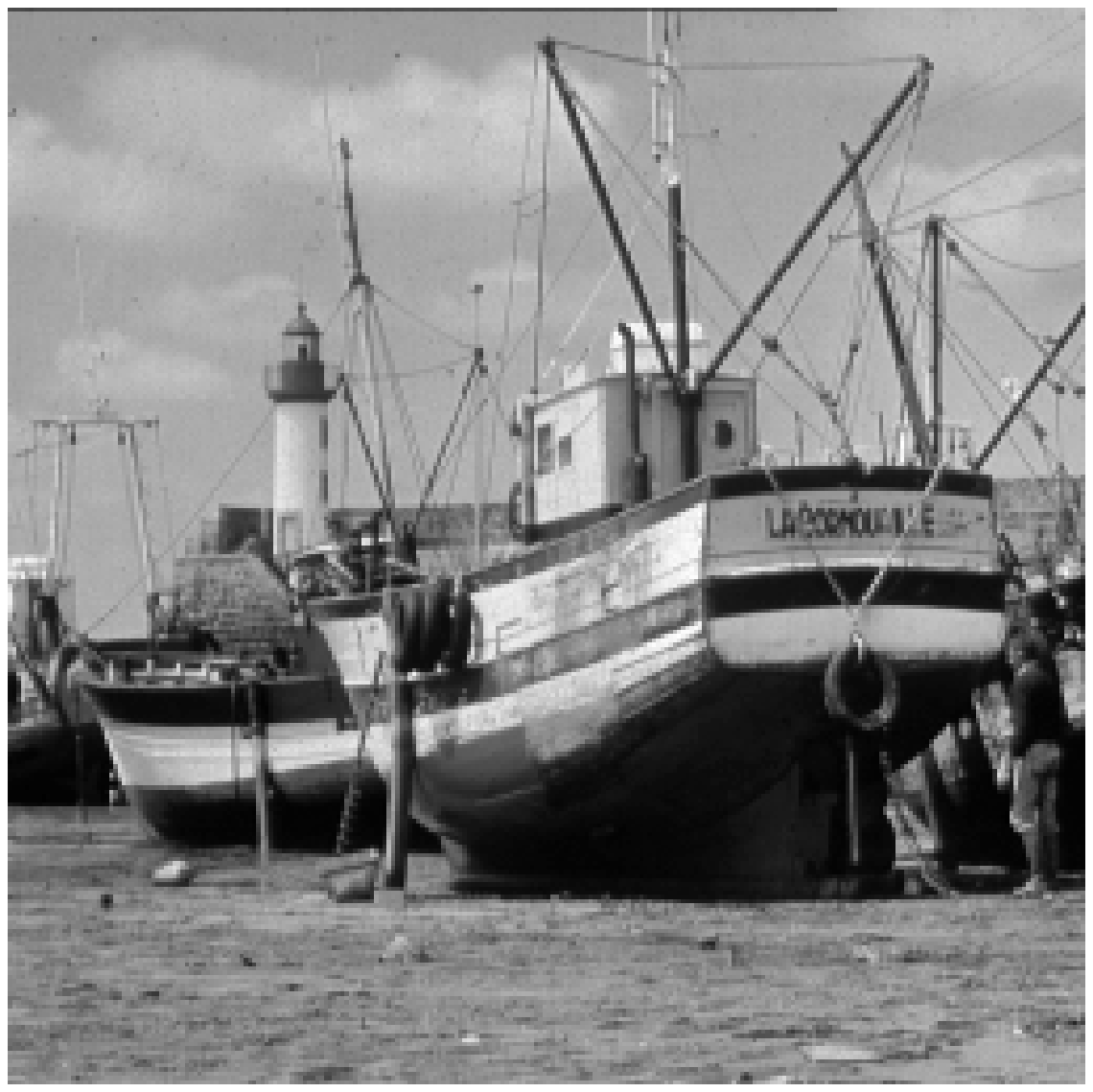}
{(e) Barchettas}
\end{minipage}\\
\caption{\label{fig1}Test images in numerical implementations.}
\end{figure*}

We use the ``Lena'' image of different sizes, i.e., $128\times128$, $256\times256$, $512\times512$ and another four images of size $256\times256$ in the numerical experiments, which are shown in Figure \ref{fig1}. To standardize the discussions, we first normalize the pixel values of the test image $\bar{f}$ to [0,255] by using the linear-stretch formula as $f=255\times\left. (\bar{f}-min(\bar{f}))\middle/(max(\bar{f})-min(\bar{f}))\right.$, where $max$ and $min$ represent the maximum and minimum of $\bar{f}$, respectively.

For the sake of simplicity, we use $\sigma$ to denote standard deviation of the white Gaussian noise and $G(hsize,\sigma)$  denotes the symmetric Gaussian low-pass filter of size $hsize$ with standard deviation $\sigma$.  In order to compare the visual perception and quality metric point of view, the performance of each method is evaluated in terms of  signal to noise ratio ($SNR$) and structural similarity index ($SSIM$): the higher $SNR$ and $SSIM$ the better the restoration results. In addition, we explicitly give the update scheme of the TMM \cite{4} and the FPM \cite{11,8} as follows
\begin{eqnarray}
&&\frac{u^{k+1}-u^k}{dt}=\left(\lambda K^T(Ku^k-f)-E_{\alpha}(u^k)\right)\hspace{10pt}(\rightarrow TMM)\\
&&(\lambda K^TK+E_{\alpha}(u^k))u^{k+1}=\lambda K^Tf \hspace{45pt}(\rightarrow  FPM)
\end{eqnarray}
by choosing a suitable original value $u^0$.  
We set $\alpha=10^{-2}$ for all of numerical implementations. Note the matrix operator in the left of the FPM is symmetric and positive definite. Therefore, we employ the conjugate gradient method to solve it as \cite{11,8}.

We first analyze the performance of our PDM by comparing it with the TMM and FPM for restoration of Lena images with different sizes. Here, the blur and noisy images are corrupted by additive white Gaussian noise with $\sigma=10$ and the Gaussian blur with $G(21,0.6)$. Table \ref{tab21} illustrates the values of $SNR$, $SSIM$ and CPU time when the numerical algorithms stop. We can observe that $SNR$ and CPU time increases while $SSIM$ decreases when the image size increases. Besides, large parameter $\lambda$ is required to penalize the data fitting term when we increase the size of the test image without increasing the level of noises. From the comparisons of $SNR$, $SSIM$, $Time$ and $Ite$, our PDM is shown better than both the TMM and FPM, especially the CPU time. In fact, the TMM needs more iterations to obtain the steady solution and the FPM needs to solve the linear equation by the numerical methods such as the conjugated gradient method with inner iterations. All of these mean that our proposed PDM is more suitable to deal with the large scale image than the other two methods.
\begin{table}[htbp]
\centering
\footnotesize{\begin{tabular}{|c|c|c|c||c|c|c||c|c|c|}
\hline {Image}&\multicolumn{9}{c|} {Lena image contaminated by  noise with $\sigma=10$.} \\
\hline {(S,$\lambda$)}&\multicolumn{3}{c||} {$(128\times128,0.14)$}&\multicolumn{3}{c||}{$(256\times256,0.12)$} &\multicolumn{3}{c|}{$(512\times512,0.11)$} \\
\hline{Method}&TMM&FPM&PDM&TMM&FPM&PDM&TMM&FPM&PDM\\
\hline SNR&17.8059&17.7039&\textbf{18.1492} &18.8162&18.7139&\textbf{19.2730}&19.9109&19.8497&\textbf{20.3225}\\
\hline SSIM&0.8229&0.8397&\textbf{0.8421}            &0.7271&0.7272&\textbf{0.7464}              &0.6230&0.5907&\textbf{0.6277}\\
\hline Time(s)&1.6224&4.3056&\textbf{0.6864} &8.1433&23.7746&\textbf{3.0108}  &68.6716&126.0488&\textbf{14.8981}\\
\hline Ite&59&\textbf{16}&32                  &79&\textbf{18}&36                              &138&\textbf{19}&38\\
\hline
\hline  {Image}&\multicolumn{9}{c|} {Lena image contaminated by  noise with $\sigma=10$  and  Gaussian blur with $G(21.0.6).$}\\
\hline {(S,$\lambda$)}&\multicolumn{3}{c||} {$(128\times128, 0.3)$}&\multicolumn{3}{c||}{$(256\times256,0.25)$} &\multicolumn{3}{c|}{$(512\times512, 0.19)$} \\
\hline {Method}&TMM&FPM&PDM&TMM&FPM&PDM&TMM&FPM&PDM\\
\hline {SNR} &14.8131&15.3534&\textbf{15.3787}&15.6258&\textbf{16.8751} &16.2984&18.2138&18.7663&\textbf{18.7726}\\
\hline {SSIM}&0.7568 &0.7919 &\textbf{0.7902} &0.6160 &\textbf{0.6903}  &0.6514  &0.5595&0.5663&\textbf{0.5778}\\
\hline {Time(s)}&1.7784 &23.6654&\textbf{2.6520} &4.6020 &101.1822&\textbf{2.6988}  &56.6596&403.1534&\textbf{19.9057}\\
\hline {Ite}&51&\textbf{29}&91&31&\textbf{30}&21&102&\textbf{24}&37\\
\hline
\end{tabular}}
\caption{\label{tab21}The related data by restoring contaminated ``Lena'' image of different sizes.  Here, $\lambda$ and  "Ite" denote the regularization parameter and iteration number, respectively, and "Time" denotes the CPU time described by second. }
\end{table}

Next, we test our PDM on other degraded images, where each test image has its own specialty, i.e., ``House" has much more sharp edges,  ``Cameraman" has more affine regions, and ``Peppers" is a relatively smooth image. To generate the degraded images, we add the additive white Gaussian noise and apply the Guassian convolution to the images. Similar results are obtained on these three test images, as shown in Table \ref{tab22}. It is observed that our PDM is more efficient than the other two methods.

\begin{table}[htbp]
\centering
\footnotesize{\begin{tabular}{|c|c|c|c||c|c|c||c|c|c|}
\hline {Description }&\multicolumn{9}{c|} {Restore noisy images generated by the second row of Figure \ref{fig1}} \\
\hline {Image}&\multicolumn{3}{c||}{House } &\multicolumn{3}{c||} {Cameraman }&\multicolumn{3}{c|}{Peppers } \\
\hline {$(\sigma,\lambda)$}&\multicolumn{3}{c||} {$(15,0.06)$}&\multicolumn{3}{c||}{$(18,0.06)$} &\multicolumn{3}{c|}{$(20,0.06)$} \\
\hline{Method}&TMM&FPM&PDM&TMM&FPM&PDM&TMM&FPM&PDM\\
\hline SNR&18.0276&17.5888&\textbf{18.1996}    &16.4785&16.5229&\textbf{17.1333} &15.4266&15.7958&\textbf{15.7982}\\
\hline SSIM&\textbf{0.4280}&0.3641&0.4117      &0.4268&0.4091&\textbf{0.4383}    &0.6083&0.6440&\textbf{0.6333}\\
\hline Time(s)&52.5255&42.4635&\textbf{6.0060} &50.0919&42.5103&\textbf{5.8812} &40.8723&38.7350&\textbf{5.6160}\\
\hline Ite&500&\textbf{21}&66                  &500&\textbf{22}&68 &415&\textbf{21}&66\\
\hline
\hline  {Description}&\multicolumn{9}{c|} {Restore noisy and blur images generated by the second row of Figure \ref{fig1}}\\
\hline {Image}&\multicolumn{3}{c||}{House }&\multicolumn{3}{c||} {Cameraman  }&\multicolumn{3}{c|}{Peppers} \\
\hline {$(\sigma,G,\lambda)$}&\multicolumn{3}{c||} {$(15,G(21,0.6),0.07)$}&\multicolumn{3}{c||}{$(18,G(21,0.6),0.08)$} &\multicolumn{3}{c|}{$(20,G(21,0.6),0.07)$} \\
\hline {Method}&TMM&FPM&PDM&TMM&FPM&PDM&TMM&FPM&PDM\\
\hline {SNR} &\textbf{16.9536}&16.5206&16.7436  &14.4791& 14.2598& \textbf{14.5334} &18.2138&18.7663&\textbf{18.7726}\\
\hline {SSIM}&\textbf{0.3831}&0.3309&0.3624     &0.3732&0.3597&\textbf{0.3760}  &0.5595&0.5663&\textbf{0.5778}\\
\hline {Time(s)}&50.9811&93.7410&\textbf{7.1136} &43.8831&97.4850&\textbf{5.7876}  &56.6596&403.1534&\textbf{19.9057}\\
\hline {Ite}&461&\textbf{26}&66   &389&\textbf{26}&54 &102&\textbf{24}&37\\
\hline
\end{tabular}}
\caption{\label{tab22}The related data by restoring contaminated image ``House'', ``Cameraman'' and ``Peppers''. }
\end{table}

Finally, we test our PDM with different smoothing parameters $\alpha$ to validate the effect of $\alpha$ in (\ref{23}). We use the image ``Barchetta'', Figure \ref{fig1} (e), and generate the degraded images by the white Gaussian noise with $\sigma=10$ and the Gaussian blurring with $G(21,0.6)$.  As shown in Table \ref{tab23}, we can obtain the overall best numerical results when $\alpha=0.01$. Furthermore, it is worthy to point out that both the TM and FPM can not be used when the model (\ref{23}) degenerates to the classic ROF model. We can use the PDM to solve (\ref{23}) as did in \cite{7} when $\alpha=0$ since the PDM does not depend on the smoothing of the model. By testing  our PDM with different $\alpha$, we find that $\alpha=0.01$ is the best choice, which gives comparable or better results than the ROF model.
\begin{table}[htbp]
\centering
\footnotesize{\begin{tabular}{|c|c||c|c|c||c|c|c||c|c|c|}
\hline { Description}&\multicolumn{10}{c|} {Barchetta image contaminated by  noise with $\sigma=10$} \\
\hline {$(\alpha,\lambda)$}&\multicolumn{1}{c||}{(0,0.19) }&\multicolumn{3}{c||}{(0.001,0.18) } &\multicolumn{3}{c||} {(0.01,0.16)}&\multicolumn{3}{c|}{(0.1,0.20) } \\
\hline{Method}&PDM&TMM&FPM&PDM&TMM&FPM&PDM&TMM&FPM&PDM\\
\hline SNR&18.6239&17.9886&18.5400&\textbf{18.6234}    &18.0055&18.3863&18.6264 &18.0459&18.6002&18.6216\\
\hline SSIM&0.7537&0.7309&0.7540&\textbf{0.7540}      &0.7352&0.7503&0.7537   &0.7340&0.7543&0.7537\\
\hline Time(s)&3.1563&29.2344&28.7031&\textbf{3.1719} &20&19.7031&2.4688&8.2188&14.4688&63.6406\\
\hline Ite&39&232&17&\textbf{34}                 &135&17&28             &52&15&500\\
\hline
\hline  {Description}&\multicolumn{10}{c|} {Barchetta image contaminated by  noise with $\sigma=10$ and Gaussian blurring with $G(21,0.6)$.}\\
\hline {$(\alpha,\lambda)$}&\multicolumn{1}{c||}{(0,0.33) }&\multicolumn{3}{c||}{(0.001,0.33) } &\multicolumn{3}{c||} {(0.01,0.33)}&\multicolumn{3}{c|}{(0.1,0.31) } \\
\hline {Method}&PDM&TMM&FPM&PDM&TMM&FPM&PDM&TMM&FPM&PDM\\
\hline {SNR} &15.7181&14.3768&15.7041&\textbf{15.7188}       &15.2739&15.7061&15.7186 &15.3533&15.7340&15.7039\\
\hline {SSIM}&0.6722&0.6126&0.6728&\textbf{0.6723}       &0.6520&0.6730&0.6715  &0.6588&0.6730&0.6732\\
\hline {Time(s)}&11.5938&6.4688&108.1719&\textbf{11.9063} &15.6250&67.4688&8.9219  &19.0313&9.1719&58.6406\\
\hline {Ite}&103&28&24&\textbf{104}   &86&23&75       &117&6&500\\
\hline
\end{tabular}}
\caption{\label{tab23}The related data by restoring contaminated ``Barchettas'' images with different values of $\alpha$.  }
\end{table}

\section{Conclusions}

We presented an image restored model based on the minimized surface regularization, which closely relates to the smoothing ROF model \cite{4}.  By using the property  of conjugate function, we first reformulate the proposed model as a min-max problem and use the primal-dual method \cite{7} to solve the optimization problem.  Theoretical convexity conditions guarantee the proposed algorithm converges to a unique global minimizer. Numerical experiments demonstrate that the proposed method holds the potential for efficient and stable computation by compared to the classic time marching method (TMM) \cite{4} and the lagged diffusivity fixed point method (FPM) \cite{11,8}, especially for the large-scale image. In the future, we would like to extend the proposed method to other image processing problem such as image inpainting, reconstruction, registration and also for vector value images, etc.

 \section*{Acknowledgements}
The authors acknowledge the financial support by the NSF of China (Nos.U1304601,11401170), Foundation of Henan Educational Committee of China (No.14A110018) and  the National Basic Research Program of China (973 Program)(No.2015CB856003).


\bibliographystyle{elsarticle-num}

\end{document}